\input amstex
%\input xy
%\input epsf
%\xyoption{all}
\documentstyle{amsppt}
\document
\magnification=1200
\NoBlackBoxes
\nologo
\hoffset1cm
\voffset2cm
\pageheight {16cm}

\def\em{\it}

%\hfill{\it Version 14}

\centerline{\bf KOLMOGOROV COMPLEXITY}

\smallskip

\centerline{\bf AS A HIDDEN FACTOR OF SCIENTIFIC DISCOURSE:}

\smallskip

\centerline{\bf FROM NEWTON'S LAW TO DATA MINING\footnotemark1} 
\footnotetext{Talk at the Plenary Session of the Pontifical Academy of Sciences on
``Complexity and Analogy in Science: Theoretical, Methodological and Epistemological Aspects'' , Casina Pio IV,
November 5--7, 2012.}
\smallskip

\bigskip

\centerline{\bf Yuri I. Manin}

\medskip

\centerline{\it Max--Planck--Institut f\"ur Mathematik, Bonn, Germany,}

\bigskip

\centerline{\bf Summary}

\medskip

 The word "complexity" is most often used as a meta--linguistic expression 
  referring to  certain intuitive characteristics of a natural
system and/or its scientific description. These characteristics may include:
sheer amount of data that must be taken into account; visible ``chaotic'' character
of these data and/or space distribution/time evolution of a system etc.
\smallskip
This talk is centered around the precise mathematical notion of  "Kolmogorov complexity",
originated in the early theoretical computer science and
measuring {\it the degree to which an available information can be compressed.}
\smallskip
 In the first part, I will argue that  a characteristic feature of basic scientific theories, from Ptolemy's
 epicycles to the Standard Model of elementary particles, is
 their splitting into two very distinct parts: the part of relatively small
 Kolmogorov complexity ("laws", "basic equations", "periodic table", "natural selection, genotypes, mutations")
 and another part, of indefinitely large Kolmogorov complexity ("initial and
 boundary conditions", "phenotypes", ``populations''). The data constituting this
 latter part are obtained by planned observations, focussed experiments, and afterwards collected
 in growing databases (formerly known as "books", "tables", "encyclopaedias" etc).
 In this discussion Kolomogorov complexity plays a role of the central metaphor.
 \smallskip
 The second part and Appendix 1 are dedicated to more precise definitions
 and examples of complexity.

 \smallskip

Finally, the last part briefly touches upon attempts to deal directly with Kolmogorov
complex massifs of data and  the ``End of Science'' prophecies.

\newpage

\centerline{\bf 1. Bi--partite structure of scientific theories}

\medskip

 In this section, I will understand the notion of  ``compression of information'' 
 intuitively and illustrate its pervasive character with several examples from the history
 of science.
 
 \smallskip
 
{\bf Planetary movements.} Firstly,  I will briefly remind the structure of several models of planetary motions 
 in the chronological order of their development.

\smallskip

After the discovery that among the stars observable by naked eye on a night sky there exist
several exceptional "moving stars" (planets), there were proposed several successful models
of their movement  that allowed predict the future positions of the moving stars.

\smallskip

The simplest of them placed all fixed stars on one celestial sphere that rotated 
around the earth in a way reflecting nightly and annual visible motions. The planets,
according to Apollonius of Perga (3rd century B.~C.), Hipparchus of Rhodes, and
 Ptolemy of Alexandria (2nd century A.~D.), were moving  in a more complicated way:
 along circular ``epicyles'' whose centers moved  along another system of circles,
 ``eccentrics''  around Earth. Data about radii of eccentrics and epicycles
 and the speed of movements were extracted from observations of the visible movements, and the
 whole model was then used in order to predict the future positions at any given
 moment of observation.  
 
 \smallskip
 
 As  D.~Park remarks ([Pa], p. 72),  ``[...]  in the midst  of all this empiricism sat the ghost of Plato, legislating that
 the curves drawn must be circles and nothing else, and that the planets
 and the various connectiong points must move along them uniformly and in no other way.''
 
 \smallskip
 
 Since in reality observable movements of planets involved accelerations, backward movements, etc.,
 two circles in place of one for each planet at least temporarily saved face of philosophy.
 Paradoxically, however, much later and much more developed mathematics of modernity
 returned to the image of ``epicycles'', that could since then form an arbitrarily high hierarchy: the idea
 of Fourier series and, later, Fourier integral transformation
 does exactly that!

\smallskip

It is well known, at least in general outline, how Copernicus replaced these geocentric models by a
heliocentric one, and how with the advent of Newton's 
$$
gravity\ law\ F=G\dfrac{m_1m_2}{r^2},\quad 
dynamic\ law\ F=ma,
$$
and the resulting solution of the ``two--body problem'', planets ``started moving''
along ellipsoidal orbits (with Sun as one focus rather than center). It is less well known to the general public that
this approximation as well is valid only insofar as we can consider negligible the gravitational forces
with which the planets interact among themselves.

\smallskip

If we  intend to obtain a more precise picture, we have to consider the system of differential
equations defining the set of curves parametrized by time $t$ in the {\it $6n$--dimensional phase space}
where $n$ is the number of planets (including Sun) taken in consideration:
$$
\frac{d^2q_i}{dt^2}=\sum_{i=1}^{n} \frac{m_im_j((q_i-q_j)}{|q_i-q_j|^3}
$$
Both Newton laws are encoded in this system.

\smallskip
The choice of one curve, corresponding to the evolution of our Solar system, 
is made when we input  {\it initial conditions} $q_i(0),  \dfrac{dq_i}{dt} (0)$ 
at certain moment of time $t=0$; they are supplied, with a certain precision, by observations.

\smallskip

At this level, a new complication emerges. Generic solutions of this system of equations, in the case of three and more bodies,
cannot be expressed by any  simple formulas (unlike the equations themselves). Moreover,
even qualitative behavior of solutions depends in extremely sensitive way on the initial conditions:
very close initial positions/velocities may produce widely divergent trajectories.
Thus, the question whether our Solar system will persist for the next, say, $10^8$ years
(even without disastrous external interventions) cannot be solved unless we know its
current parameters (masses of planets, positions of their centers of mass, and speeds) with unachievable precision. 
This holds even without appealing to much more precise Einstein's description of gravity,
or without taking in account comets, asteroid belts and Moons of the Solar system
(the secondary planets turning around planets themselves).

\smallskip

It goes without saying that a similarly detailed description of, say, our Galaxy, 
taking in account movements of all individual
celestial bodies, constituting it, is unachievable from the start,
because of sheer amount of  these bodies. Hence, to understand its general space--time structure, we must first construct models
involving averaging on a very large scale. And of course, the model of space--time itself,
now involving Einstein's equations, will describe an ``averaged'' space--time.

\medskip

{\bf  Information compression: first summary.} In this brief summary of consecutive scientific models,
one can already see the following persisting pattern: the subdivision into a highly compressed part (``laws'')
and potentionally indefinitely complex part.
The first part in our brief survey was represented by formulas that literally
became cultural symbols of Western civilization: Newton's laws, that were followed by 
Einstein's $E=mc^2$ and Heisenberg's $pq-qp=\dfrac{h}{2\pi i}.$ The second part is kinematically represented
by ``initial'' or ``boundary'' conditions, and dynamically by a potentially
unstable character of dependence of the data we are interested in from these initial/boundary conditions.

\smallskip

More precisely, a mathematical description of the ``scene'' upon which develops kinematics and dynamics in these models is also represented by highly compressed
mathematical images, only this time of geometric nature. Thus, the postulate that kinematics of a single massive point
is represented by its position in an ideal Euclidean space represents one of the "laws" as well.
To describe kinematics, one should amplify this ``configuration space''  and replace it
by the ``phase space'' parametrizing positions and velocities, or, better, momenta. For 
one massive point it is  a space of dimension {\it six}:
this is  the answer of mathematics to  Zeno's  ``Achilles and the Turtle'' paradox.
For a planet system consisting of $n$ planets (including Sun) the phase space has
dimension $6n$.

\smallskip

For Einstein's equations of gravitation, the relevant picture is much more complicated: it involves
configuration and phase spaces that have {\it infinite} dimension,
and require quite a fair amount of mathematics for their exact description.
Nevertheless, this part of our models is still clearly separated from the one that we refer to
as the part of infinite Kolmogorov complexity, because mathematics developed
a concise language for description of geometry.

\smallskip

One more lesson of  our analysis is this: ``laws'' can be discovered and efficiently used only
if and when we restrict our attention to definite domains, space--time scales, and
kinds of matter and interactions.  For example, there was no place for chemistry in the pictures above.

\medskip

{\bf From macroworld to microworld: the Standard Model of elementary particles
and interactions.} From astronomy, we pass now to the deepest known level of microworld:
theory of elementary particles and their interactions.

\smallskip

I will  say a few words about the so called Standard Model of the elementary particles and their interactions, that
took its initial form in the 1970's as a theoretical construction in the framework of the Quantum Field Theory.
The Standard Model got its first important experimental correlates with the discovery of
quarks (components of nuclear ``elementary'' particles) and $W$ and $Z$ bosons, quanta of ineractions.
For a very rich and complex history of this stage of theoretical physics, stressing the role of experiments and experimenters,
see the fascinating account [Zi] by Antonio Zichichi. The Standard Model recently reappeared
on the first pages of the world press thanks to the renewed hopes that the  last critically
missing  component of the Model, the Higgs boson, has  finally been observed.

\smallskip

Somewhat paradoxically, one can say that mathematics of the Standard Model is firmly based 
on the same ancient archetypes of
the human  thought as  that of Hipparchus and Ptolemy: symmetry and uniform
movement along circles.

\smallskip

More precisely, the basic idea of symmetry of modern classical (as opposed to quantum) non--relativistic physics   
involves the symmetry group of rigid movements of the three--dimensional Euclidean space,
that is combinations of parallel shifts and rotations around a point. The
group of rotations is denoted $SO(3)$, and celestial spheres are the unique 
objects invariant with respect to rotations. Passing from Hipparchus and Ptolemy to
modernity includes two decisive steps:  adding shifts (Earth, and then Sun, cease being centers of the Universe),
and, crucially,  understanding the new meta--law of physics:  {\it symmetry} must govern {\it laws of physics themselves}
rather than objects/processes etc that these laws are supposed to govern (such as Solar System).

\smallskip

When we pass now to the quantum mechanics, and further to the Quantum Field Theory (not involving gravitation),
the group of $SO(3)$ (together with shifts) should be extended, in particular, by several copies of such groups as 
$SU(2)$ and $SU(3)$ describing rotations in the {\it internal degrees of freedom} of elementary
particles, such as spin, colour etc. The basic "law" that should be invariant
with respect to this big group, is encoded in the Lagrangian density: it is a ``mathematical formula''
that is considerably longer than everything we get exposed to in our high school
and even college courses: cf. Appendix 2.

\smallskip

Finally, the Ptolemy celestial movements, superpositions of rotations of rigid spheres, now 
transcends our space--time and happens in the infinite--dimensional Hilbert space
of wave--functions: this is the image describing, say, a hydrogen atom in the paradigm
of the first decades of the XXth century.

\bigskip

{\bf Information compression: second summary.}  I will use the examples above in order to justify the following
viewpoint.

\smallskip

Scientific laws  (at least those that are expressed by mathematical constructions) can be considered as
{\it programs for computation}, whereas observations produce {\it inputs} to these programs.

\smallskip
 {\it Outputs}  of these computations serve first to check/establish a domain of applicability of our theories.
We compare the predicted behavior of a system with observed one, we are happy when our predictions
agree quantitatively and/or qualitatively with observable behaviour, we fix the border signs
signalling that at this point we went too far.

\smallskip
 Afterwards, the outputs
are used for practical/theoretical purposes, e.~g. in engineering,  weather predictions etc,
but also to formulate the new challenges arising before the scientific thinking.
\smallskip

This comparison of scientific laws with programs is, of course, only a metaphor, but
it will allow us to construct also a precise model of the kind of complexity,  inherently associated
with this metaphor of science: Kolmogorov complexity.

\smallskip

The next section is dedicated to the sketch of  this notion in the framework of
mathematics, again in its historical perspective.

\bigskip
\centerline{\bf 2.  Integers and their Kolmogorov complexity}

\medskip

{\bf Positional notations as programs.} In this section, I will explain that the well known to the general public decimal notations
of natural numbers are themselves programs. 

\smallskip

What are they supposed to calculate?

\smallskip

Well, the actual {\it numbers}   that are encoded by this notation, and are more adequately represented
by, say, rows of strokes:

$$
7 :\quad |||||||,\quad\quad 13:\quad ||||||||||||| ,\quad ...\quad\quad ,\, 1984:\quad ||||...||||
$$
Of course, in the last example it is unrealistic even to expect that if I type here 1984 strokes,
an unsophisticated reader will be able to check that I am not mistaken. There will be simply too much strokes
to count, whereas the {\it notation--program} ``1984'' contains only four signs chosen from the alphabet of ten signs.
One can save on the size of alphabet, passing to the binary notation, then ``1984'' will be replaced by a longer program
``11111000000''.  However, comparing the length of the program with the ``size'' of the number, i.~e. the
respective  number of strokes, we see that decimal/binary notation gives an immense economy:
the program length is approximately the logarithm of the number of strokes (in the base 10 or 2 respectively).

\smallskip

More generally, we can speak about ``size'', or ``volume'' of any finite text based upon a fixed finite alphabet.

\smallskip

The discovery of this logarithmic upper bound of the Kolmogorov complexity of {\it numbers} was
a leap in the development of humanity on the scale of civilizations.

\smallskip

However, if one makes some slight additional conventions in the system of notation,
it will turn out that {\it some} integers admit a much shorter notation.  For example,
let us allow ourselves to use the vertical dimension and write, e.~g.  $10^{10^{10}}$.

\smallskip

The logarithm of the last number is about $10^{10}$, much larger than the length of the notation for which we used only 6 signs!
And if we are unhappy about non--linear notation, we may add to the basic alphabet two brackets (,) and
postulate that $a(b)$ means $a^b$. Then   $10^{10^{10}}$ will be linearly written as $10(10(10))$ using only 10 signs,
still much less than $10^{10}+1$ decimal digits (of course, $10^{10}$ of them will be just zeroes).

\smallskip

Then, perhaps, {\it all} integers can be produced by notation/programs that are much shorter than
logarithm of their size?

\smallskip

No! It turns out that absolute majority of numbers (or texts) {\it cannot} be significantly compressed, although an infinity of integers
can be written in a much shorter way than it can be done in any chosen system of  positional notation.

\smallskip

If we leave the domain of integers and leap, to, say, such a number as $\pi =3,1415926...$, it looks as if it had infinite complexity.
However, this is not so. There exists a program that can take as input the (variable) place of a decimal digit
(an integer) and give as output the respective digit. Such a program is itself a text in a chosen algorithmic language, and as such, it
also has a complexity: its own Kolmogorov complexity.  One agrees that this is the complexity of $\pi$.
\smallskip

A reader should be aware that I have left many subtle points of the definition of Kolmogorov complexity in shadow,
in particular, the fact that its dependence of the chosen system of encoding and computation model
can change it only by a bounded quantity etc.
A reader  who would like to see some more mathematics about this matter is referred to the brief Appendix 1
and the relevant references.

\smallskip

Here I will mention two other remarkable facts related to the Kolmogorov complexity of numbers: one regarding its
unexpected relation to the idea of {\it randomness}, and another one showing that some psychological data
make explicit the role of  this complexity in the cognitive activity of our mind.

\medskip

{\bf Complexity and randomness.} Consider arbitrarily long finite sequences of zeroes and ones,
say, starting with one so that each such sequence could be interpreted as a binary notation
of an integer. 

\smallskip

There is an intuitive notion of ``randomness '' of such a sequence. In the contemporary technology
``random'' sequences of digits and similar random objects are used for encoding information, in order to make
it inaccessible for third parties. In fact, a small distribured industry producing such random sequences (and, say, random big primes)
has been created. A standard way to produce random objects is to leave mathematics and
to recur to physics: from throwing a piece to registering white noise.

\smallskip

One remarkable property of Kolmogorov complexity is this: {\it those sequences of digits whose Kolmogorov
complexity is approximately the same as their length, are random in any meaningful sense of the word.} In particular,
they cannot be generated by a program essentialy shorter  than the sequence itself.

\medskip

{\bf Complexity and human mind.} In the history of humanity, discovery  of  laws of 
classical and quantum physics   that represent incredible compression of complex information,
stresses the role of  Kolmogorov complexity, at least as a relevant metaphor for understanding the laws
of cognition.

\smallskip

In his very informative book [De], Stanislas Dehaene considers certain experimental results about the
statistics of appearance numerals and other names of numbers. cf. especially pp. 110 -- 115,
subsection ``Why  are some numerals more frequent than others?''.

\smallskip

As mathematicians, let us consider the following abstract question: can one say anything non-obvious about
possible probabilities distributions on the set of {\it all} natural numbers? More precisely,
one such distribution is a sequence of non--negative real numbers $p_n, n=1,2,\dots $ such that
$\sum_n p_n =1.$  Of course, from the last formula it follows that $p_n$ must tend to zero, when $n$
tends to infinity; moreover $p_n$ cannot tend to zero too slowly: for example, $p_n=n^{-1}$ will not do.
But two different distributions can be widely incomparable.

\smallskip

Remarkably, it turns out that if we restrict our class of distributions only to {\it computable from below} ones,
that is, those in which $p_n$ can be computed as a function of $n$ (in a certain precise sense),
then it turns out that there is a distinguished and small subclass $C$ of such distributions,
that are in a sense  {\it maximal} ones. Any member $(p_n)$ of this class has the following unexpected property
(see [Lev]):

\smallskip

{\it the probability $p_n$ of the number $n$, up to a bounded (from above and below) factor, equals the inverse
of the exponentiated Kolmogorov complexity of $n$.}

 \smallskip

This statement needs additional qualifications: the most important one is that we need here
{\it not} the original Kolmogorov complexity but the so called {\it prefix--free} version of it.
We omit technical details, because they are not essential here.  But the following
properties of any distribution $(p_n)\in C$ are worth stressing in our context:

\smallskip

(i) Most of the numbers $n$, those that are Kolmogorov "maximally complex",appear with
probability comparable with $n^{-1}\,({\roman{log}\,n})^{-1-\varepsilon}$, with a small $\varepsilon$:
``most large numbers appear with frequency inverse to their size'' (in fact, somewhat smaller one).

\smallskip

(ii) However, frequencies of  those numbers that are Kolmogorov very simple, such as $10^3$ (thousand),
$10^6$ (million), $10^9$ (billion), produce sharp local peaks in the graph of $(p_n)$.

\smallskip

The reader may compare these properties of the  discussed class of distributons, which can be
called {\it a priori distributions}, with the observed frequencies of numerals (number words)
in printed and oral texts in various languages: cf. Dehaene, loc. cit., p.~111,  Figure 4.4.
To me, their qualitative agreement looks very convincing: brains and their societies
do reproduce a priori probabilities. 

\smallskip

Notice that  those parts of the Dehaene and Mehler graphs in loc.~cit. that refer
to large numbers, are somewhat misleading: they might create an impression that frequencies of the numerals, say,
between $10^6$ and $10^9$ smoothly interpolate between those of   $10^6$ and $10^9$ themselves,
whereas in fact they abruptly drop down.

\smallskip

Finally, I want to stress that the class of a priori probability distributions that we are considering here
is {\it qualitatively distinct} from those that form now  a common stock of sociological and sometimes scientific analysis:
cf. a beautiful synopsis by Terence Tao in [Ta]. The appeal to the uncomputable degree of maximal compression is
exactly what can make such a distribution an eye--opener. As I have written at the end of [Ma2]:

\smallskip

``One can argue that all cognitive activity of our
civilization, based upon symbolic (in particular, mathematical)
representations of reality, deals actually with the {\it initial
Kolmogorov segments}
of potentially infinite linguistic constructions,
{\it always} replacing  vast volumes of data by their
compressed descriptions. This is especially visible
in the outputs of the modern genome projects.

\smallskip

In this sense, such linguistic cognitive activity
can be metaphorically compared to a gigantic precomputation 
process, shellsorting infinite worlds of expressions
in their Kolmogorov order.''

\bigskip

\centerline{\bf 3. New cognitive toolkits: WWW and databases}

\medskip

{\bf ``The  End of Theory''.} In summer 2008, an issue of the ``Wired Magazine'' appeared.
It's cover story ran: ``The End of Theory: The Data Deluge Makes the Scientific Method
Obsolete''.

\medskip

The message of this essay, written by the Editor--in--Chief Chris Anderson,
was summarized in the following words:

\smallskip

{\it ``The new availablility of huge amounts of data, along with statistical tools
to crunch these numbers, offers a whole new way of understanding the world.
Correlation supersedes causation, and science can advance even without coherent 
models, unified theories, or really any mechanical explanation at all. There's 
no reason to cling to our old ways. It's time to ask: What can science learn from Google?"}

\medskip

I will return to this rhetoric question at the end of this talk.
Right now I want only to stress that, as well as in the scientific models of
the ``bygone  days'', basic theory is unavoidable  in this brave new Petabyte World:
encoding and decoding data, search algorithms, and of course, computers themselves
are just engineering embodiment of some very basic and very abstract
notions of mathematics. The mathematical idea underlying 
the structure of modern computers is the Turing machine (or one of several other
equivalent formulations of the concepts of computability). We know that the
universal Turing machine has a very small Kolmogorov complexity,
and therefore, using the basic metaphor of this talk, we can say that 
the bipartite structure of the classical scientific theories is reproduced
at this historical stage.
 
\smallskip

Moreover, what Chris Anderson calls ``the new availability of huge amounts of data''
by itself is not very new: after spreading of printing, astronomic
observatories, scientific laboratories, and statistical studies, the amount of data available
to any visitor of a big public library  was always huge, and studies of correlations
proliferated for at least the last two centuries.

\smallskip
Charles Darwin himself collected the database of his observations, and the result
of his pondering over it was the theory of evolution.
\smallskip
A representative recent example is the book  [FlFoHaSCH],
sensibly reviewed in [Gr].

\smallskip

Even if the sheer volume of data has by now grown
by several orders of magnitude, this is not the gist of Anderson's
rhetoric.
\smallskip
What Anderson actually wants to say is that human beings are now -- happily! -- free from thinking
over these data. Allegedly, computers will take this burden upon themselves, and will
provide us with correlations -- replacing the old--fashioned ``causations'' (that I prefer to call
scientific laws) -- and expert guidance.

\smallskip
Leaving aside such questions as  how ``correlations'' might possibly help us  understand the structure of Universe or
predict the Higgs boson, I would like to quote the precautionary tale from [Gr]:

\smallskip

{\it   ``[...] in 2000 Peter C.~Austin, a medical statistician at the University
of Toronto, and his colleagues conducted a study of all 10,674,945 residents
of Ontario aged between eighteen and one hundred. Residents were
randomly assigned to different groups, in which they were classified according to their
astrological signs. The research team then searched through more than two hundred of the most
common diagnoses of hospitalization until they identified two where patients 
under one astrological sign had a significantly higher probability of hospitalization
compared to those born under the remaining signs combined: Leos had a higher
probability of gastrointestinal hemorrage while Sagittarians had a higher
probability of fracture of the upper arm compared to all other signs combined.

\smallskip

It is thus relatively easy to generate statistically significant but spurious correlations 
when examining a very large data set and a similarly large number of potential variables.
Of course, there is no biological mechanism whereby Leos might be predisposed to intestinal
bleeding or Sagittarians to bone fracture, but Austin notes,   `It is tempting to construct
biologically plausible reasons for observed subgroup effects after having observed them.'
Such an exercise is termed `data mining', and Austin warns,
`Our study therefore serves as a cautionary note regarding the interpretation of findings
generated by data mining' [...]''}

\bigskip

\centerline{\bf Coda}

\medskip

{\bf What can science learn from Google:}

\smallskip
\centerline{ ``Think! Otherwise no Google will help you.''}

\newpage

\centerline{\bf References}

\medskip
[An] Ch.~Anderson. {\it The End of Theory.} in: Wired, 17.06, 2008.

\smallskip

[ChCoMa] A.~H.~Chamseddine, A.~Connes, M.~Marcolli. {\it Gravity and the standard model with
neutrino mixing.} Adv.Theor.Math.Phys. 11 (2007), 991-1089. arXiv:hep-th/0610241
\smallskip

[De]  S.~Dehaene. {\it The Number Sense. How the Mind creates Mathematics.}
Oxford UP, 1997.
\smallskip
[FlFoHaSCH] R.~Floud, R.~W.~Fogel, B.~Harris, Sok Chul Hong. {\it The  Changing Body:
Health, Nutrition and Human Development in the
Western World Since 1700.} Cambridge UP, 431 pp.
\smallskip

[Gr] J.~Groopman. {\it The Body and the Human Progress.}
In: NYRB, Oct. 27, 2011
\smallskip
[Lev] L.~A.~Levin, {\em Various measures of complexity for finite objects (axiomatic
description)}, Soviet Math. Dokl. Vol.17 (1976) N. 2, 522--526.

\smallskip
[Ma1] Yu.~Manin. {\it A Course of Mathematical Logic for Mathematicians.}
2nd Edition, with new Chapters written by Yu.~Manin and B.~Zilber. Springer, 2010.
\smallskip

[Ma2] Yu.~Manin. {\it Renormalization and Computation II: Time Cut-off and the Halting Problem.}
 {\it In: Math. Struct. in Comp. Science,} pp. 1--23, 2012, Cambridge UP.
Preprint math.QA/0908.3430
\smallskip
[Pa] D.~Park. {\it The How and the Why. An Essay on the Origins and Development of Physical Theory.} 
Princeton UP, 1988.
\smallskip

[Ta] T.~Tao. {\it E pluribus unum: From Complexity, Universality.} Daedalus, Journ. of the AAAS, Summer 2012, 
23--34.

\smallskip
[Zi] A.~Zichichi. {\it Subnuclear Physics. The first 50 years: Highlights from Erice to ELN}.
World Scientific, 1999.

\newpage
\centerline{{\it APPENDIX 1.} {\bf A brief guide to computability}}

\medskip

This appendix contains a sketch of the mathematical computability theory, or theory of algorithmic
computations, as it was born in the first half the XXth century in the work of such thinkers as
Alonso Church, Alan Turing and Andrei Kolmogorov. We are not concerned here with its
applied aspects studied under the general heading of ``Computer Science''.

\smallskip

This theory taught us two striking lessons. 
\smallskip
 First, that there is a unique universal notion 
of computability in the sense that all seemingly very different versions of it turned out to be equivalent.
We will sketch here the form that is called {\it the theory
of (partial) recursive functions.}

\smallskip
Second, that this theory sets its own limits and unavoidably leads to confrontation with {\it uncomputable
problems.}

\smallskip

Both discoveries led to very interesting research aiming to the extension of this territory
of classical computability such as basics of theory of quantum computing etc. But we are not concerned
with this development here.

\smallskip
{\bf Three descriptions of partial recursive functions}. The subject of the
theory of recursive functions is a set of functions whose domain and values are vectors of natural numbers
of  arbitrary fixed lengths:   $f:\bold{Z}_+^m\to \bold{Z}_+^n$.

\smallskip
 
An important qualification: a ``function'',  say,  $f:\,X\to Y$, below always means a pair   $(f, D(f))$, where
$D(f)\subset X$ and $f$ is a map of sets $D(f)\to Y$. The definition domain $D(f)$ is not always mentioned
explicitly. If $D(f)=X$, the function might  be called ``total''; generally it may be called ``partial'' one.

\medskip

{\bf (i) Intuitive description.} A function $f:\bold{Z}_+^m\to \bold{Z}_+^n$
is called {\it (partial) recursive} iff it is ``semi--computable'' in the following sense:

\medskip
 there exists an algorithm $F$ accepting as inputs vectors $x=(x_1,\dots ,x_m)\in \bold{Z}_+^m$.
 with the following properties:
 
 \smallskip
 
 -- if $x\in D(f)$,  $F$ produces as output $f(x)$.
 
 \smallskip
 
 -- if $x\notin D(f)$, $F$ either produces answer ``NO'', or works indefinitely long without
producing any output.

\medskip

{\bf (ii) Formal\ description\ (sketch).} It starts with two lists:

\smallskip

-- An explicit list of ``obviously'' semi--computable {\it basic
functions} such as constant functions,
projections onto $i$--th coordinate etc.

\smallskip

-- An explicit list of {\it elementary\ operations}, such as an inductive definition, that can be applied to several
semi--computable functions and ``obviously'' produces from them
a new semi--computable function. 

\smallskip
The key elementary operation involves finding
the least root of equation $f(x)=y$ (if it exists) where $f$ is already defined function.
It is this operation that involves search and makes introduction of partial functions inevitable.

\smallskip

-- After that, the set of {\it partial  recursive  functions}  is defined as {\it the minimal set of functions
containing all basic functions and closed wrt all elementary operations.}

\medskip

{\bf (iii) Diophantine description} (A DIFFICULT THEOREM).   A function $f:\bold{Z}_+^m\to \bold{Z}_+^n$
is partial recursive iff there is a polynomial 
$$
P(x_1,\dots, x_m;y_1,\dots ,y_n; t_1,\dots ,t_q)\in \bold{Z}[x,y,t]
$$
such that the graph 
$$
\Gamma_f := \{(x,f(x))\} \subset  \bold{Z}_+^m\times \bold{Z}_+^n
$$
  is the projection  of the subset  $P=0$ in
$\bold{Z}_+^m\times \bold{Z}_+^n\times \bold{Z}_+^q$.

\medskip

{\bf Constructive worlds.} {\it An (infinite) constructive world} is a countable set  $X$
(usually of some finite Bourbaki structures) given together with a class of  
{\it structural\ numberings:}
intuitively computable bijections
$\nu :\bold{Z}_+\to X$ which form a principal homogeneous space over the group
of  recursive permutations of  $\bold{Z}_+$. An element $x\in X$
is called {\it a constructive object.}

\smallskip

{\it Example:} $X$ = all finite words in a fixed finite alphabet $A$.

\bigskip

{\bf Church's thesis:} {\it  Let $X$, $Y$ be two constructive worlds,
$\nu_X :\bold{Z}_+\to X$,  $\nu_Y :\bold{Z}_+\to Y$ their structural numberings, and  $F$ an (intuitive) algorithm 
that takes as input an object  $x\in X$ and produces an
object $F(x)\in Y$ whenever $x$ lies in the domain of definition of $F$.

\smallskip

Then  $f:= \nu_Y^{-1}\circ F\circ\nu_X:\bold{Z}_+\to \bold{Z}_+$ is a partial recursive function.}

\medskip

The status of Church's thesis in mathematics is very unusual. It is not a theorem, since it involves
an undefined notion of  ``intuitive computability''. It expresses the fact that several developed
mathematical constructions  whose explicit goal was to formalize this notion, led to provably equivalent results.
But moreover, it expresses the belief that any new such attempt will inevitably produce
again an equivalent notion. Briefly, Church's thesis is ``an experimental fact in the mental world''.

\medskip

{\bf Exponential Kolmogorov complexity of constructive objects.}
 Let $X$ be a constructive world. For any (semi)--computable function $u:\,\bold{Z}_+\to X$,
the (exponential)  complexity of an object $x\in X$ relative to $u$ is
$$
K_u(x):= \roman{min}\,\{m\in  \bold{Z}_+\,|\,u(m)=x\}.
$$
If such $m$ does not exist, we put   $K_u(x)= \infty .$

\medskip
 {\bf Claim:} there exists such  $u$ (``an optimal Kolmogorov numbering'', or ``decompressor'') that
for each other $v$, some constant $c_{u,v}>0$, and all $x\in X$,
$$
K_u(x)\le c_{u,v}K_v(x).
$$

 This $K_u(x)$ is called  {\it exponential Kolmogorov\ complexity} of $x$.
\medskip

{\it A Kolmogorov order} on a constructive world
$X$ is a bijection  $\bold{K}=\bold{K}_u:\, X\to \bold{Z}_+$ arranging elements of $X$
in the increasing order of their complexities $K_u$.

\medskip

The reader must keep in mind two warnings related to these notions:

\medskip

-- Any optimal numbering is only partial function, and its definition domain is not decidable.

\medskip

-- Kolmogorov complexity $K_u$ itself is {\it not\ computable.} It is the lower bound
of a sequence of computable functions. Kolmogorov order is not computable as well.

\medskip

-- Kolmogorov order of $\bold{Z}_+$ {\it cardinally differs  from  the  natural order}
 in the following sense:
it puts in the initial segments very large numbers that are at the same time Kolmogorov simple.

\smallskip

 {\it Example}: let $a_n:= n^{n^{.^{.^{.^{n}}}}}$ ($n$ times).  
\smallskip
Then $K_u(a_n)\le cn$ for some $c>0$ and all $n>0$.

\medskip

{\bf Kolmogorov complexity of recursive functions.} When we spoke about ``complexity of $\pi$'',
we had in mind a program that, given an input $n\in \bold{Z}_+$, calculates the $n$--th decimal
digit of $\pi$. Such a program calculates therefore a recursive function. However, the set of all recursive functions
$f:\bold{Z}_+^m\to \bold{Z}_+^n$ {\it do not} form a constructive world, if $m>0$!

\smallskip

Nevertheless, one can speak about a Kolmogorov optimal enumeration of programs, computing functions of this set,
and  thus about Kolmogorov complexity of such functions themselves.  This is critically
important for applicability of our  ``complexity metaphor'' in the domain of scientific knowledge.
In fact, both "laws" and descriptions of ``phase spaces'' that make the scene for these laws
belong rather to domains of intuitively computable functions than to the domain of constructive objects.

\smallskip

Mathematical details of constructions, underlying brief explanations collected in this Appendix can be found in [Ma1],
Chapter II.

\newpage
\centerline{{\it APPENDIX 2.} {\bf Lagrangian of the Standard Model}}

\centerline{{\it Source:  [ChCoMa]}}

\medskip

``In flat space and in Lorentzian signature the Lagrangian of the
standard model with neutrino mixing and Majorana mass terms, written
using the Feynman gauge fixing, is of the form

$$
{\Cal{L}}_{SM}=
$$
$$-\frac{1}{2}\partial_{\nu}g^{a}_{\mu}\partial_{\nu}g^{a}_{\mu}
-g_{s}f^{abc}\partial_{\mu}g^{a}_{\nu}g^{b}_{\mu}g^{c}_{\nu}
-\frac{1}{4}g^{2}_{s}f^{abc}f^{ade}g^{b}_{\mu}g^{c}_{\nu}g^{d}_{\mu}g^{e}_{\nu}
-\partial_{\nu}W^{+}_{\mu}\partial_{\nu}W^{-}_{\mu}-M^{2}W^{+}_{\mu}W^{-}_{\mu}-
$$
$$
\frac{1}{2}\partial_{\nu}Z^{0}_{\mu}\partial_{\nu}Z^{0}_{\mu}-\frac{1}{2c^{2}_{w}}
M^{2}Z^{0}_{\mu}Z^{0}_{\mu}-\frac{1}{2}\partial_{\mu}A_{\nu}\partial_{\mu}A_{\nu}
-igc_{w}(\partial_{\nu}Z^{0}_{\mu}(W^{+}_{\mu}W^{-}_{\nu}-
-W^{+}_{\nu}W^{-}_{\mu})
-Z^{0}_{\nu}(W^{+}_{\mu}\partial_{\nu}W^{-}_{\mu}-
$$
$$
W^{-}_{\mu}\partial_{\nu}W^{+}_{\mu})
+Z^{0}_{\mu}(W^{+}_{\nu}\partial_{\nu}W^{-}_{\mu}-W^{-}_{\nu}\partial_{\nu}W^{+}_{\mu}))
-igs_{w}(\partial_{\nu}A_{\mu}(W^{+}_{\mu}W^{-}_{\nu}-W^{+}_{\nu}W^{-}_{\mu})
-A_{\nu}(W^{+}_{\mu}\partial_{\nu}W^{-}_{\mu}-
$$
$$
W^{-}_{\mu}\partial_{\nu}W^{+}_{\mu})
+A_{\mu}(W^{+}_{\nu}\partial_{\nu}W^{-}_{\mu}-W^{-}_{\nu}\partial_{\nu}W^{+}_{\mu}))
-\frac{1}{2}g^{2}W^{+}_{\mu}W^{-}_{\mu}W^{+}_{\nu}W^{-}_{\nu}+\frac{1}{2}g^{2}
W^{+}_{\mu}W^{-}_{\nu}W^{+}_{\mu}W^{-}_{\nu}+
$$
$$
g^2c^{2}_{w}(Z^{0}_{\mu}W^{+}_{\mu}Z^{0}_{\nu}W^{-}_{\nu}-Z^{0}_{\mu}Z^{0}_{\mu}W^{+}_{\nu}
W^{-}_{\nu})
+g^2s^{2}_{w}(A_{\mu}W^{+}_{\mu}A_{\nu}W^{-}_{\nu}-A_{\mu}A_{\mu}W^{+}_{\nu}
W^{-}_{\nu})
+g^{2}s_{w}c_{w}(A_{\mu}Z^{0}_{\nu}(W^{+}_{\mu}W^{-}_{\nu}-
$$
$$
W^{+}_{\nu}W^{-}_{\mu})
-2A_{\mu}Z^{0}_{\mu}W^{+}_{\nu}W^{-}_{\nu})
%% higgs kinetic terms from minimal coupling with gauge
-\frac{1}{2}\partial_{\mu}H\partial_{\mu}H-2M^2\alpha_{h}H^{2}
-\partial_{\mu}\phi^{+}\partial_{\mu}\phi^{-}
-\frac{1}{2}\partial_{\mu}\phi^{0}\partial_{\mu}\phi^{0} -
$$
$$
\beta_{h}\left(\frac{2M^{2}}{g^{2}}+\frac{2M}{g}H+\frac{1}{2}(H^{2}+\phi^{0}\phi^{0}+2\phi^{+}\phi^{-
})\right)
  +\frac{2M^{4}}{g^{2}}\alpha_{h}
  %% higgs quartic potential
  -g\alpha_h
M\left(H^3+H\phi^{0}\phi^{0}+2H\phi^{+}\phi^{-}\right)-
$$
$$
\frac{1}{8}g^{2}\alpha_{h}
\left(H^4+(\phi^{0})^{4}+4(\phi^{+}\phi^{-})^{2}
+4(\phi^{0})^{2}\phi^{+}\phi^{-}
+4H^{2}\phi^{+}\phi^{-}+2(\phi^{0})^{2}H^{2}\right)
%% minimal coupling higgs-gauge L_gH
-gMW^{+}_{\mu}W^{-}_{\mu}H-
$$
$$
\frac{1}{2}g\frac{M}{c^{2}_{w}}Z^{0}_{\mu}Z^{0}_{\mu}H
-\frac{1}{2}ig\left(W^{+}_{\mu}(\phi^{0}\partial_{\mu}\phi^{-}
-\phi^{-}\partial_{\mu}\phi^{0})
-W^{-}_{\mu}(\phi^{0}\partial_{\mu}\phi^{+}
-\phi^{+}\partial_{\mu}\phi^{0})\right)+
$$
$$
\frac{1}{2}g\left(W^{+}_{\mu}(H\partial_{\mu}\phi^{-}
-\phi^{-}\partial_{\mu}H)
 +W^{-}_{\mu}(H\partial_{\mu}\phi^{+}-\phi^{+}\partial_{\mu}H)\right)
+\frac{1}{2}g\frac{1}{c_{w}}(Z^{0}_{\mu}(H\partial_{\mu}\phi^{0}-\phi^{0}\partial_{\mu}H)+
$$
$$
M\,(\frac{1}{c_{w}}Z^{0}_{\mu}\partial_{\mu}\phi^{0}+W^{+}_{\mu}
\partial_{\mu}\phi^{-}+W^{-}_{\mu}
\partial_{\mu}\phi^{+})
-ig\frac{s^{2}_{w}}{c_{w}}MZ^{0}_{\mu}(W^{+}_{\mu}\phi^{-}-W^{-}_{\mu}\phi^{+})
   +igs_{w}MA_{\mu}(W^{+}_{\mu}\phi^{-}-W^{-}_{\mu}\phi^{+})-
$$
$$
ig\frac{1-2c^{2}_{w}}{2c_{w}}Z^{0}_{\mu}(\phi^{+}\partial_{\mu}\phi^{-}
-\phi^{-}\partial_{\mu}\phi^{+})
+igs_{w}A_{\mu}(\phi^{+}\partial_{\mu}\phi^{-}-\phi^{-}\partial_{\mu}\phi^{+})
-\frac{1}{4}g^{2}W^{+}_{\mu}W^{-}_{\mu}
\left(H^{2}+(\phi^{0})^{2}+2\phi^{+}\phi^{-}\right) -
$$
$$
\frac{1}{8}
g^{2}\frac{1}{c^{2}_{w}}Z^{0}_{\mu}Z^{0}_{\mu}
\left(H^{2}+(\phi^{0})^{2}+2(2s^{2}_{w}-1)^{2}\phi^{+}\phi^{-}\right)
-\frac{1}{2}g^{2}\frac{s^{2}_{w}}{c_{w}}Z^{0}_{\mu}\phi^{0}(W^{+}_{\mu}\phi^{-}+W^{-}_{\mu}\phi^{+})
$$
$$
-\frac{1}{2}ig^{2}\frac{s^{2}_{w}}{c_{w}}Z^{0}_{\mu}H(W^{+}_{\mu}\phi^{-}-W^{-}_{\mu}\phi^{+})
+\frac{1}{2}g^{2}s_{w}A_{\mu}\phi^{0}(W^{+}_{\mu}\phi^{-}+W^{-}_{\mu}\phi^{+})
+\frac{1}{2}ig^{2}s_{w}A_{\mu}H(W^{+}_{\mu}\phi^{-}-W^{-}_{\mu}\phi^{+})-
$$
$$
g^{2}\frac{s_{w}}{c_{w}}(2c^{2}_{w}-1)Z^{0}_{\mu}A_{\mu}\phi^{+}\phi^{-}
-g^{2}s^{2}_{w}A_{\mu}A_{\mu}\phi^{+}\phi^{-}
%% kinetic terms for fermions L_gf
+\frac 12 i
g_s\,\lambda_{ij}^a(\bar{q}^{\sigma}_{i}\gamma^{\mu}q^{\sigma}_{j})g^{a}_{\mu}
-\bar{e}^{\lambda}(\gamma\partial+m^{\lambda}_{e})e^{\lambda}
-\bar{\nu}^{\lambda}(\gamma\partial+
$$
$$
m^{\lambda}_{\nu})\nu^{\lambda}
-\bar{u}^{\lambda}_{j}(\gamma\partial+
m^{\lambda}_{u})u^{\lambda}_{j}
-\bar{d}^{\lambda}_{j}(\gamma\partial+m^{\lambda}_{d})d^{\lambda}_{j}
+igs_{w}A_{\mu}\left(-(\bar{e}^{\lambda}\gamma^{\mu}
e^{\lambda})+\frac{2}{3}(\bar{u}^{\lambda}_{j}\gamma^{\mu}
u^{\lambda}_{j})-\frac{1}{3}(\bar{d}^{\lambda}_{j}\gamma^{\mu}
d^{\lambda}_{j})\right)+
$$
$$
\frac{ig}{4c_{w}}Z^{0}_{\mu}
\{(\bar{\nu}^{\lambda}\gamma^{\mu}(1+\gamma^{5})\nu^{\lambda})+
(\bar{e}^{\lambda}\gamma^{\mu}(4s^{2}_{w}-1-\gamma^{5})e^{\lambda})
     +(\bar{d}^{\lambda}_{j}\gamma^{\mu}(\frac{4}{3}s^{2}_{w}-1-\gamma^{5})d^{\lambda}_{j})+
(\bar{u}^{\lambda}_{j}\gamma^{\mu}(1-\frac{8}{3}s^{2}_{w}+\gamma^{5})u^{\lambda}_{j})
\}+
$$
$$
+\frac{ig}{2\sqrt{2}}W^{+}_{\mu}\left((\bar{\nu}^{\lambda}\gamma^{\mu}(1+\gamma^{5})U^{lep}_{\lambda\kappa}e^{\kappa})
+(\bar{u}^{\lambda}_{j}\gamma^{\mu}(1+\gamma^{5})C_{\lambda\kappa}d^{\kappa}_{j})\right)
$$
$$
+\frac{ig}{2\sqrt{2}}W^{-}_{\mu}\left((\bar{e}^{\kappa}U^{lep\dagger}_{\kappa\lambda}\gamma^{\mu}(1+\gamma^{5})\nu^{\lambda})
+(\bar{d}^{\kappa}_{j}C^{\dagger}_{\kappa\lambda}\gamma^{\mu}(1+\gamma^{5})u^{\lambda}_{j})\right)+
$$
$$
   %% yukawa coupling of higgs with fermions (modulo some mass terms) L_Hf
   +\frac{ig}{2M\sqrt{2}}\phi^{+}
\left(-m^{\kappa}_{e}(\bar{\nu}^{\lambda}U^{lep}_{\lambda\kappa}(1-\gamma^{5})e^{\kappa})
+m^{\lambda}_{\nu}(\bar{\nu}^{\lambda}U^{lep}_{\lambda\kappa}(1+\gamma^{5})e^{\kappa}\right)+
$$
$$
   +\frac{ig}{2M\sqrt{2}}\phi^{-}
\left(m^{\lambda}_{e}(\bar{e}^{\lambda}U^{lep\dagger}_{\lambda\kappa}(1+\gamma^{5})\nu^{\kappa})
-m^{\kappa}_{\nu}(\bar{e}^{\lambda}U^{lep\dagger}_{\lambda\kappa}(1-\gamma^{5})\nu^{\kappa}\right)
-\frac{g}{2}\frac{m^{\lambda}_{\nu}}{M}H(\bar{\nu}^{\lambda}\nu^{\lambda})
-\frac{g}{2}\frac{m^{\lambda}_{e}}{M}H(\bar{e}^{\lambda}e^{\lambda})+
$$
$$
\frac{ig}{2}\frac{m^{\lambda}_{\nu}}{M}\phi^{0}(\bar{\nu}^{\lambda}\gamma^{5}\nu^{\lambda})
-\frac{ig}{2}\frac{m^{\lambda}_{e}}{M}\phi^{0}(\bar{e}^{\lambda}\gamma^{5}e^{\lambda})
-\frac 14\,\bar
\nu_\lambda\,M^R_{\lambda\kappa}\,(1-\gamma_5)\hat\nu_\kappa -\frac
14\,\overline{\bar
\nu_\lambda\,M^R_{\lambda\kappa}\,(1-\gamma_5)\hat\nu_\kappa}+
$$
$$
\frac{ig}{2M\sqrt{2}}\phi^{+}
\left(-m^{\kappa}_{d}(\bar{u}^{\lambda}_{j}C_{\lambda\kappa}(1-\gamma^{5})d^{\kappa}_{j})
+m^{\lambda}_{u}(\bar{u}^{\lambda}_{j}C_{\lambda\kappa}(1+\gamma^{5})d^{\kappa}_{j}\right)+
$$
$$
 \frac{ig}{2M\sqrt{2}}\phi^{-}
\left(m^{\lambda}_{d}(\bar{d}^{\lambda}_{j}C^{\dagger}_{\lambda\kappa}(1+\gamma^{5})u^{\kappa}_{j})
-m^{\kappa}_{u}(\bar{d}^{\lambda}_{j}C^{\dagger}_{\lambda\kappa}(1-\gamma^{5})u^{\kappa}_{j}\right)
-\frac{g}{2}\frac{m^{\lambda}_{u}}{M}H(\bar{u}^{\lambda}_{j}u^{\lambda}_{j})
-\frac{g}{2}\frac{m^{\lambda}_{d}}{M}H(\bar{d}^{\lambda}_{j}d^{\lambda}_{j})+
$$
$$
\frac{ig}{2}\frac{m^{\lambda}_{u}}{M}\phi^{0}(\bar{u}^{\lambda}_{j}\gamma^{5}u^{\lambda}_{j})
-\frac{ig}{2}\frac{m^{\lambda}_{d}}{M}\phi^{0}(\bar{d}^{\lambda}_{j}\gamma^{5}d^{\lambda}_{j})
$$
\bigskip

Here the notation is  as follows:

\medskip
$\bullet$  Gauge bosons: $A_\mu, W_\mu^{\pm}, Z_\mu^0, g^{a}_{\mu}$

$\bullet$ Quarks: $u^{\kappa}_{j}, d^{\kappa}_{j}$, collective : $ q^{\sigma}_{j}$

$\bullet$ Leptons: $e^{\lambda}, \nu^{\lambda}$

$\bullet$ Higgs fields: $H, \phi^{0}, \phi^{+}, \phi^{-}$

$\bullet$ Ghosts: $G^a, X^{0}, X^{+}, X^{-}, Y$,

$\bullet$ Masses: $m^{\lambda}_{d}, m^{\lambda}_{u}, m^{\lambda}_{e}, m_h,
M$ (the latter is the mass of the $W$)

$\bullet$ Coupling constants $g=\sqrt{4\pi \alpha}$ (fine structure), $g_s=$
strong, $\alpha_h=\frac{m_h^2}{4M^2}$

$\bullet$ Tadpole Constant $\beta_h$

$\bullet$ Cosine and sine of the weak mixing angle $c_w, s_w$

$\bullet$ Cabibbo--Kobayashi--Maskawa mixing matrix: $C_{\lambda\kappa}$

$\bullet$ {Cabibbo-Kobayashi-Maskawa matrix}

$\bullet$  Structure constants of ${SU}(3)$: $f^{abc}$

$\bullet$ The Gauge is the Feynman gauge.''

\enddocument

\bigskip
\centerline{\bf References}

\medskip
[An] Ch.~Anderson. {\it The End of Theory.} in: Wired, 17.06, 2008.
\smallskip
[ChCoMa] A.~H.~Chamseddine, A.~Connes, M.~Marcolli. {\it Gravity and the standard model with
neutrino mixing.} Adv.Theor.Math.Phys. 11 (2007), 991-1089. arXiv:hep-th/0610241
\smallskip
[De]  S.~Dehaene. {\it The Number Sense. How the Mind creates Mathematics.}
Oxford UP, 1997.
\smallskip

[FlFoHaSCH] R.~Floud, R.~W.~Fogel, B.~Harris, Sok Chul Hong. {\it The  Changing Body:
Health, Nutrition and Human Development in the
Western World Since 1700.} Cambridge UP, 431 pp.
\smallskip

[Gr] J.~Groopman. {\it The Body and the Human Progress.}
In: NYRB, Oct. 27, 2011.

\smallskip

[Lev] L.~A.~Levin, {\em Various measures of complexity for finite objects (axiomatic
description)}, Soviet Math. Dokl. Vol.17 (1976) N. 2, 522--526.

\smallskip
[Ma] Yu.~Manin. {\it A Course of Mathematical Logic for Mathematicians.}
2nd Edition, with new Chapters written by Yu.~Manin and B.~Zilber. Springer, 2010.

\smallskip
[Pa] D.~Park. {\it The How and the Why. An Essay on the Origins and Development of Physical Theory.} 
Princeton UP, 1988.
\smallskip
[Zi] A.~Zichichi. {\it Subnuclear Physics. The first 50 years: Highlights from Erice to ELN}.
World Scientific, 1999.

\enddocument